\documentclass[11pt]{article}
\usepackage{amsmath,amssymb,comment}
  \usepackage{paralist}
  \usepackage{graphics} 
  \usepackage{epsfig} 
\usepackage{graphicx}  \usepackage{epstopdf}
 \usepackage[colorlinks=true]{hyperref}
\hypersetup{urlcolor=blue, citecolor=red}

\begin{document}


\begin{center}

{\LARGE\textbf{Benjamin-Ono model of an internal wave under a flat surface \\}} \vspace {10mm}
\vspace{1mm} \noindent

{\large \bf Alan C. Compelli $^{a}$} 
\\
{\it E-mail address:  alan.compelli@ucc.ie} 

\hskip 1cm 

{\large \bf Rossen I. Ivanov $^{b}$} 
\\
{\it E-mail address: rossen.ivanov@dit.ie} 

\hskip 1cm 

\vskip 1cm

$^a${\it School of Mathematical Sciences, University College Cork,  }{\it Cork, Ireland} 
\\

$^b${\it School of Mathematical Sciences, Dublin Institute of Technology, }{\it Kevin street, Dublin 8, Ireland} 
\hskip 1cm
\\

\vskip 1cm
\end{center}

\begin{abstract}A two-layer fluid system separated by a pycnocline in the form of an internal wave is considered. The lower layer is infinitely deep, with a higher density than the upper layer which is bounded above by a flat surface. The fluids are incompressible and inviscid. A Hamiltonian formulation for the fluid dynamics is presented and it is shown that an appropriate scaling leads to the integrable Benjamin-Ono equation.
\end{abstract}

\section{Introduction} Due to the ubiquitous presence of waves on the surface of oceans much research has been completed which attempts to explain the mechanisms responsible for such waves. However, despite many centuries of reporting, by mariners in particular, of unexplained phenomena beneath the surface, such as the ``dead water'' phenomenon coined by Fridtjof Nansen \cite{Nansen}, it is relatively recently that some progress of a mathematical nature has been achieved in the description of internal waves. In many situations both the surface and the internal waves are moving in the presence of currents. Due to the nonlinearity the wave-current interactions are quite complex and necessitate the study of rotational fluids \cite{CompelliIvanov1JNMP, Constantin_2011, Jonsson, Peregrine, ThomasKlopman} . 

The significant findings of Vladimir Zakharov in 1968 in \cite{Zakharov} have established the use of the Hamiltonian approach for dynamic descriptions of wave motion. For single layer systems both irrotational \cite{BenjOlv, CraigGroves1,Milder, Miles} and rotational \cite{ConstantinEscher, ConstantinEscher2, NearlyHamiltonian, ConstantinSattingerStrauss, ConstantinStrauss, TelesdaSilvaPeregrine, Wahlen_HamFormWaterWaves} set-ups have been examined within the Hamiltonian framework. Stratified systems which contain a pycnocline separating two layers in the form of an internal wave have been considered in an irrotational context in \cite{BenjBridPart1,BenjBridPart2,CraigGuyenneKalisch,CraigGuyenneSul3} and in a rotational context, which is of most interest, in \cite{Compelli1Wavemotion, Compelli2MonatshMath,ConstantinIvanov,ConstantinJohnson}. 

We aim to establish a model for a system with a flat-bed, flat surface and internal wave. Our setting captures the nonlinear dynamics as influenced by both the internal wave behaviour and the presence of a depth-varying current. In physical terms the system can be thought of as a non-mixing oceanic environment consisting of two discrete fluid bodies at different densities in the presence of a current such as the case in the Pacific ocean where the EUC (Equatorial Under Current) influences a stratified region of the ocean \cite{FedorovBrown}.

Perturbative techniques can be used to develop such models. In \cite{CompelliIvanov2} the chosen scaling regime leads to a KdV type model. We choose to adopt a similar approach, but with a scaling which leads to a Benjamin-Ono type approximation.

\section{Set-up and governing equations}

\begin{figure} [h]
\centering
\includegraphics[width=0.75\textwidth]{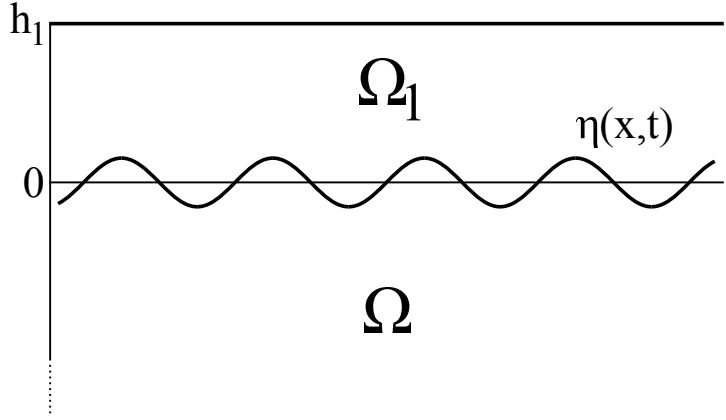}
\caption{The system under study.}
\label{fig:1}
\end{figure}

The flow we consider consists of two layers, $\Omega$ and $\Omega_1$, which have different densities, due to different salinity levels or temperatures. They are separated by a very thin layer termed as a pycnocline (whose thickness is neglected), as shown in Figure \ref{fig:1} where the horizontal axis is $x$ and the vertical axis is $y$. The internal wave is formed at the pycnocline $y=\eta(x,t)$. The mean of $\eta$ is assumed to be zero, that is
\begin{equation} \label{int_eta}
 \int\limits_{\mathbb{R}} \eta(x,t) dx=0
 \end{equation}
for all $t$ so that $\eta(x,t)$ measures the elevation of the internal wave with respect to the level $y=0$.

The domains $\Omega$ and $\Omega_1$ are defined as
\begin{align}
\Omega&:=\{(x,y)\in \mathbb{R}^2:-\infty<y<\eta(x,t)\}\notag\\
\text{and} \quad \Omega_1&:=\{(x,y)\in \mathbb{R}^2:\eta(x,t)<y< h_1\}.\notag
\end{align} 
Throughout the paper we will use the subscript 1 to denote quantities pertaining to the upper layer, while quantities referring to the lower layer will appear without subscript. The domain $\Omega$ is infinitely deep and $\Omega_1$ has a flat surface. The physical justification for this assumption of absence of surface motion is due to the small amplitude of the surface waves in comparison to the amplitude of the internal waves. The fluids are considered to be inviscid and incompressible with $\rho$ and $\rho_1$ being the respective constant densities of the lower and upper media. Stable stratification requires 
\begin{equation}
\label{stability}
\rho>\rho_1.
\end{equation}

The velocity fields are given by
\begin{equation}
{\bf{V}}(x,y,t)=(u,v)\mbox{ and }{\bf{V}}_1(x,y,t)=(u_1,v_1)
\end{equation} and the incompressibility implies $\text{div}{\bf{V}}=0 $ and $\text{div} {\bf{V}}_1 =0$  or \begin{equation} \label{div}
u_{x}+v_{y}=0  \mbox{ and } u_{1,x}+v_{1,y}=0.
\end{equation}
The Euler equations 
\begin{eqnarray}
\frac{\partial}{\partial t}{\bf{V}}+({\bf{V}}\cdot \nabla){\bf{V}} &=& -\frac{1}{\rho}\nabla p + {\bf g},  \label{Eu1}\\
\frac{\partial}{\partial t}{\bf{V}}_1+({\bf{V}}_1\cdot \nabla){\bf{V}}_1 &=& -\frac{1}{\rho_1}\nabla p_1 + {\bf g} \label{Eu2}
\end{eqnarray} govern the fluid dynamics in each layer, where $p, p_1$ are the dynamic pressure terms in the corresponding layers and ${\bf g}=(0, -g)$ is the Earth's gravitational acceleration.

In addition, the kinematic boundary condition for the interface between the layers (the pycnocline) must be satisfied, that is
\begin{equation} \label{BC}
 \eta_t =v-\eta_x u = v_1 - \eta_x u_{1} \quad \text{on} \quad y=\eta(x,t).
 \end{equation}

Velocity potentials $\varphi$ and $\varphi_1$ are introduced to capture the irrotational components of the velocity fields. These components in general include time-independent constant currents so that the velocity potentials can be further decomposed as 
\begin{equation}
\varphi \equiv \widetilde{\varphi}+\kappa x\mbox{ and }\varphi_1 \equiv\widetilde{\varphi}_1+\kappa_1 x\notag
\end{equation}
where $\kappa$ and $\kappa_1$ are the respective time-independent currents at $y=0$. In this decomposition $\widetilde{\varphi}$ and $\widetilde{\varphi}_1$
represent the potential components related to the wave motion. In the presence of a current with constant vorticity we can therefore represent the velocity field components in $\Omega$ via
\begin{equation}
 u  =   \widetilde{\varphi}_{x} +\gamma y+\kappa\mbox{ and }  v=  { \widetilde{\varphi}}_{y},
\end{equation}
and similarly for
$\Omega_1$ via
\begin{equation}
u_1  =   \widetilde{\varphi}_{1,x} +\gamma_1 y+ \kappa_1\mbox{ and }v_1=  { \widetilde{\varphi}}_{1,y}
\end{equation}
where $\gamma=u_y-v_x$ and $\gamma_1=u_{1,y}-v_{1,x}$ are the constant non-zero vorticities (see \cite{ConstantinIvanovMartin}, where the situation with a free surface is studied).

The incompressibility of the fluid \eqref{div} allows us to introduce a stream function $\psi$ in $\Omega$ via
\begin{equation}
 u  =\psi_y\mbox{ and }  v= -\psi_x
\end{equation}
and a stream function $\psi_1$ in 
$\Omega_1$ via
\begin{equation}
u_1  =\psi_{1,y}  \mbox{ and }  v_1=-\psi_{1,y}.
\end{equation}

 The following assumptions will be made: $\eta(x,t)<h_1$ for all values of $x $ and $t$;  $\eta(x,t)$ and $\widetilde{\varphi}_1(x,y,t)$, are in the Schwartz class $\mathcal{S}(\mathbb{R})$ with respect to the $x$ variable (for any $y$ and $t$); $\widetilde{\varphi}(x,y,t)$ is in the Schwartz class $\mathcal{S}(\mathbb{R}^2)$ with respect to both the $x$ and $y$ variables (for any $t$). Due to these assumptions for large absolute values of $x$ the internal wave attenuates, meaning that
 \begin{equation}
\lim_{|x|\rightarrow \infty}\eta(x,t)=\lim_{|x|\rightarrow \infty}{ \widetilde{\varphi}}(x,y,t)=\lim_{|x|\rightarrow \infty}{ \widetilde{\varphi}_1}(x,y,t)=0.
\end{equation}
Moreover, \begin{equation}
\lim_{y\rightarrow -\infty}{ \widetilde{\varphi}}(x,y,t)=0
\end{equation} for all values of $x$ and $t$. The physical reasoning for this is the absence of wave motion at infinite depth $y\to -\infty$.
 
Euler's equations \eqref{Eu1}, \eqref{Eu2} in terms of the introduced variables are \cite{CompelliIvanov2}
\begin{align}
\label{F1_GenSystem}
&\widetilde{\varphi}_{t} +\frac{1}{2}| \nabla \psi|^2-\gamma\psi_1+gy+\frac{p}{\rho}=f(t)\\
\label{F2_GenSystem}
\text{and} \quad &\widetilde{\varphi}_{1,t} +\frac{1}{2}| \nabla \psi_1|^2-\gamma_1 \psi_1+gy+\frac{p_1}{\rho_1}=f_1(t)
\end{align}
where $p$ and $p_1$ are the dynamic pressure terms, $\rho$ and $\rho_1$ are the constant densities, $f$ and $f_1$ are some so far arbitrary functions of time. Their presence is related to the fact that the velocity potentials are determined up to an additive term whose gradient is zero. For further convenience we choose
$$
\rho f(t)=\rho_1 f_1(t).$$

At the interface $y=\eta(x,t)$ (denoted by a subscript ``i'') the dynamic pressure terms are equal giving the Bernoulli equation \cite{Johnson_Book}
\begin{multline}
\rho\Big(({\widetilde{\varphi}_{t}})_{\mathrm{i}}+\frac{1}{2}|\nabla \psi|_{\mathrm{i}}^2-\gamma\chi +g\eta+f(t)\Big)\\=\rho_1\Big(( {\widetilde{\varphi}_{1,t}})_{\mathrm{i}}+\frac{1}{2}|\nabla \psi_1|_{\mathrm{i}}^2-\gamma_1 \chi_1 +g\eta+f_1(t)\Big)
\end{multline}
where $\chi$ and $\chi_1$ are the interface stream functions. Furthermore, it could be shown \cite{Compelli1Wavemotion} that $\chi=\chi_1$, and so the Bernoulli equation becomes
\begin{equation}
\rho\Big(({\widetilde{\varphi}_{t}})_{\mathrm{i}}+\frac{1}{2}|\nabla \psi|_{\mathrm{i}}^2-\gamma\chi +g\eta\Big)=\rho_1\Big(( {\widetilde{\varphi}_{1,t}})_{\mathrm{i}}+\frac{1}{2}|\nabla \psi_1|_{\mathrm{i}}^2-\gamma_1 \chi +g\eta\Big)
\end{equation}

or

\begin{equation} \label{E1}
( \rho{\widetilde{\varphi}_{t}}- \rho_1  {\widetilde{\varphi}_{1,t}})_{\mathrm{i}}=\frac{\rho_1}{2}|\nabla \psi_1|^2_{\mathrm{i}}-\frac{\rho}{2}|\nabla \psi|^2_{\mathrm{i}}+(\rho \gamma-\rho_1 \gamma_1) \chi +(\rho_1 - \rho) g\eta.
\end{equation}

This form suggests the introduction of a single variable $\rho{\widetilde{\varphi}}- \rho_1  {\widetilde{\varphi}_{1}}  $  which is going to provide one of the Hamiltonian coordinates (the {\it momentum}) in the next section. The other one (the {\it coordinate}) is the variable $\eta(x,t)$, which satisfies the kinematic boundary condition \eqref{BC}
\begin{equation} \label{E2}
 \eta_t =(\widetilde{\varphi}_{y})_{\mathrm{i}} - \eta_x \big(  (\widetilde{\varphi}_{x})_{\mathrm{i}} + \gamma \eta + \kappa\big)= (\widetilde{\varphi}_{1,y})_{\mathrm{i}} - \eta_x \big(  (\widetilde{\varphi}_{1,x})_{\mathrm{i}} + \gamma_1 \eta + \kappa_1\big).
 \end{equation}

\section{The Hamiltonian formulation}

The functional $H$, which describes the total energy of the system, can be written as the sum of the kinetic, $\mathcal{K}$, and potential energy, $\mathcal{V}$ contributions. The potential part, must be  $$V(\eta)=\lim_{h\to -\infty}\rho g\int\limits_{\mathbb{R}} \int\limits_{h}^{\eta}  y \, dy dx
+\rho_1 g\int\limits_{\mathbb{R}} \int\limits_{\eta}^{h_1}  y \, dy dx.$$ However, the potential energy is always measured from some reference value, e.g. $V(\eta=0)$ which is the potential energy of the current (without wave motion). Therefore, the relevant part of the potential energy, contributing to the wave motion is  $$ \mathcal{V}(\eta)= V(\eta)-V(0)=\rho g\int\limits_{\mathbb{R}} \int\limits_{0}^{\eta}  y \, dy dx
+\rho_1 g\int\limits_{\mathbb{R}} \int\limits_{\eta}^{0}  y \, dy dx=\frac{1}{2}(\rho-\rho_1)g\int\limits_{\mathbb{R}} \eta^2 dx
.$$

In order to determine the kinetic energy of the wave motion, from the total kinetic energy of the fluid
\begin{equation}
\frac{1}{2}\rho\int\limits_{\mathbb{R}} \int\limits_{-\infty}^{\eta}  (u^2+v^2)dy dx+\frac{1}{2}\rho_1\int\limits_{\mathbb{R}} \int\limits_{\eta}^{h_1}  (u_1^2+v_1^2)dy dx
\end{equation} one should subtract again the constant, but infinite kinetic energy of the current which is 
\begin{equation}
\frac{1}{2}\rho\int\limits_{\mathbb{R}} \int\limits_{-\infty}^{0}  (\gamma y+\kappa)^2dy dx+\frac{1}{2}\rho_1\int\limits_{\mathbb{R}} \int\limits_{0}^{h_1}  (\gamma_1 y+\kappa_1)^2dy dx.
\end{equation} 
 
In terms of the dependent variables $\eta(x,t)$, $\widetilde{\varphi}(x,t)$ and $\widetilde{\varphi}_1(x,t)$ this kinetic energy is 
\begin{multline}
\mathcal{K} (\eta,\widetilde{\varphi},\widetilde{\varphi}_1)=\frac{1}{2}\rho\int\limits_{\mathbb{R}} \int\limits_{-\infty}^{\eta}  \left(( \widetilde{\varphi}_{x} +\gamma y+\kappa)^2+(\widetilde{\varphi}_y)^2\right)dy dx  -    \frac{1}{2}\rho\int\limits_{\mathbb{R}} \int\limits_{-\infty}^{0}  (\gamma y+\kappa)^2dy dx\\
+\frac{1}{2}\rho_1\int\limits_{\mathbb{R}}\int\limits_{\eta}^{h_1}  \left(( \widetilde{\varphi}_{1,x} +\gamma_1 y+\kappa_1)^2+(\widetilde{\varphi}_{1,y})^2\right)dy dx -\frac{1}{2}\rho_1\int\limits_{\mathbb{R}} \int\limits_{0}^{h_1}  (\gamma_1 y+\kappa_1)^2dy dx\\
= \frac{1}{2}\rho\int\limits_{\mathbb{R}} \int\limits_{-\infty}^{\eta}  \left(( \widetilde{\varphi}_{x})^2 +(\widetilde{\varphi}_y)^2+2\widetilde{\varphi}_{x}(\gamma y+\kappa)\right)dy dx    \\
+\frac{1}{2}\rho_1\int\limits_{\mathbb{R}}\int\limits_{\eta}^{h_1}  \left(( \widetilde{\varphi}_{1,x})^2 +(\widetilde{\varphi}_{1,y})^2 + 2\widetilde{\varphi}_{1,x} (\gamma_1 y+\kappa_1)\right)dy dx \\
 +\frac{1}{6}(\rho\gamma^2-\rho_1\gamma_1^2)\int\limits_{\mathbb{R}} \eta^3dx+\frac{1}{2}\big(\rho\gamma\kappa-\rho_1\gamma_1\kappa_1\big)\int\limits_{\mathbb{R}} \eta^2dx.
\end{multline}

The Hamiltonian is therefore

\begin{multline}
 H(\eta,\widetilde{\varphi},\widetilde{\varphi}_1)=\mathcal{K}+\mathcal{V}
= \frac{1}{2}\rho\int\limits_{\mathbb{R}} \int\limits_{-\infty}^{\eta}  \left(( \widetilde{\varphi}_{x})^2 +(\widetilde{\varphi}_y)^2+2\widetilde{\varphi}_{x}(\gamma y+\kappa)\right)dy dx    \\
+\frac{1}{2}\rho_1\int\limits_{\mathbb{R}}\int\limits_{\eta}^{h_1}  \left(( \widetilde{\varphi}_{1,x})^2 +(\widetilde{\varphi}_{1,y})^2 + 2\widetilde{\varphi}_{1,x} (\gamma_1 y+\kappa_1)\right)dy dx \\
 +\frac{1}{6}(\rho\gamma^2-\rho_1\gamma_1^2)\int\limits_{\mathbb{R}} \eta^3dx
 +\frac{1}{2}\big((\rho\gamma\kappa-\rho_1\gamma_1\kappa_1)+(\rho-\rho_1)g\big)\int\limits_{\mathbb{R}} \eta^2dx. 
\end{multline}

The assumption that $\widetilde{\varphi}(x,y,t)$ is in the Schwartz class $\mathcal{S}(\mathbb{R}^2)$ with respect to both the $x$ and $y$ variables gives

$$ \lim_{x\to \pm \infty} \int_{-\infty}^0 (\gamma y+\kappa) \widetilde{\varphi}(x,y,t) dy =0 \qquad \text{for any $t,$ }$$  and furthermore
$$
\int\limits_{\mathbb{R}} \int\limits_{-\infty}^{\eta}  \widetilde{\varphi}_{x}(\gamma y+\kappa) dy dx    = -\int\limits_{\mathbb{R}}\widetilde{\varphi}(x, \eta,t) (\gamma \eta + \kappa)\eta_x dx.
 $$

We would like to write the Hamiltonian only in terms of the one dimensional variables pertaining to the interface. To this end the Dirichlet-Neumann operators
\begin{equation}
G(\eta)\phi=\left(\frac{\partial\widetilde{\varphi}}{\partial {\bf{n}}}\right)_{\mathrm{i}} \sqrt{1+\eta_x ^2}\mbox{ and } G_1(\eta)\phi_1=\left(\frac{\partial\widetilde{\varphi}_1}{\partial {\bf{n}}_1}\right)_{\mathrm{i}} \sqrt{1+\eta_x ^2}
\end{equation}
are introduced, where ${\bf{n}}$ and ${\bf{n}}_1$ are the unit exterior normals, $\sqrt{1+(\eta_x )^2}$ is a normalisation factor and
\begin{equation} \label{phi}
       \phi(x,t):=\widetilde{\varphi}(x,\eta(x,t),t) 
        \mbox{ and }
        \phi_1(x,t):=\widetilde{\varphi}_1(x,\eta(x,t),t)
\end{equation}
are the interface velocity potentials. 

Since usually there is no jump in the current velocity, in what follows we take $\kappa=\kappa_1.$

The Hamiltonian can therefore be written in terms of conjugate variables $\eta(x,t)$ and $\xi(x,t)$, following the procedure in \cite{CompelliIvanov2}, as
\begin{multline}
\label{Hamiltonian_conj}
H(\eta,\xi)=\frac{1}{2}\int\limits_{\mathbb{R}} \xi G(\eta) B^{-1}G_1(\eta)\xi \,dx
- \frac{1}{2}\rho\rho_1(\gamma-\gamma_1)^2\int\limits_{\mathbb{R}}   \eta\eta_x   B^{-1}\eta\eta_x  dx \\
-\gamma\int\limits_{\mathbb{R}} \xi\eta\eta_x dx
-\kappa\int\limits_{\mathbb{R}} \xi\eta_x  dx+\rho_1(\gamma-\gamma_1)\int\limits_{\mathbb{R}} \eta\eta_x B^{-1}G(\eta)\xi\,dx
+\frac{1}{6}(\rho\gamma^2-\rho_1\gamma_1^2)\int\limits_{\mathbb{R}} \eta^3dx\\
+\frac{1}{2}\big((\rho\gamma-\rho_1\gamma_1)\kappa+g(\rho-\rho_1 )\big)\int\limits_{\mathbb{R}} \eta^2dx 
\end{multline}
where 
\begin{alignat}{2} \label{xi}
\xi(x,t):=\rho\phi(x,t)-\rho_1\phi_1(x,t)
\end{alignat}
and the operator $B$, as per \cite{CraigGuyenneKalisch}, is introduced as
\begin{equation}
\label{B_definition}
B:=\rho G_1(\eta)+\rho_1 G(\eta).
\end{equation}

We point out that due to the initial assumptions on $\eta,\widetilde{\varphi}$ and $\widetilde{\varphi}_1 $ and \eqref{phi},\eqref{xi}, the Hamiltonian variables $\eta$ and $\xi$ are in $\mathcal{S}(\mathbb{R})$ with respect to the $x$ variable for all $t$.

The equations of motion \eqref{E1} -- \eqref{E2} can be presented in the following form (see \cite{Compelli1Wavemotion} for details)
\begin{equation}
\label{EOMsys}
\eta_t=\frac{\delta H}{\delta\xi}\mbox{ and }
\xi_t=-\frac{\delta H}{\delta\eta} +\Gamma  \chi
\end{equation} 
where 
\begin{alignat}{2}
\Gamma:=\rho\gamma-\rho_1\gamma_1.
\end{alignat}
Furthermore, we can write
\begin{alignat}{2}
\label{lem2}
\chi(x,t)=- \int_{-\infty}^x\eta_t (x',t)dx'=-\partial_x^{-1}\eta_t
\end{alignat}
and introducing a non-canonical Poisson bracket as \cite{Wahlen_HamFormWaterWaves}
\begin{equation}
\label{Poisson_Bracket_Non_Can}
\{F_1,F_2\}=\int\limits_{\mathbb{R}}\bigg(\frac{\delta F_1}{\delta \eta(x)}\frac{\delta F_2}{\delta \xi(x)} -\frac{\delta F_1}{\delta \xi(x)}\frac{\delta F_2}{\delta \eta(x)} \bigg)dx-\Gamma \int\limits_{\mathbb{R}}\bigg(\frac{\delta F_1}{\delta \xi(x)}\int\limits_{-\infty}^x\frac{\delta F_2}{\delta \xi(x')}dx'\bigg)dx
\end{equation}
for functionals $F_1$ and $F_2$ provided at least one of the functionals in addition satisfies $$  \int\limits_{\mathbb{R}}\frac{\delta F_k}{\delta \xi(x)}  dx=0,$$
 the equations (\ref{EOMsys}) can be written in the form 
 
 $$\xi_t=\{\xi, H\} , \qquad   \eta_t=\{\eta, H\}.$$

The canonical Hamiltonian form can be achieved under the transformation (cf. \cite{Compelli1Wavemotion,Compelli2MonatshMath,Wahlen_HamFormWaterWaves})
\begin{alignat}{2}
\label{vartrans}
\xi\rightarrow\zeta=\xi-\frac{\Gamma}{2} \int_{-\infty}^{x} \eta(x',t)\,dx',
\end{alignat}
giving the equations
\begin{equation}
\label{EOMsysCAN}
\eta_t=\frac{\delta H}{\delta\zeta}\mbox{ and }
\zeta_t=-\frac{\delta H}{\delta\eta}.
\end{equation} The condition \eqref{int_eta} ensures that $\int_{-\infty}^{x }\eta(x',t)dx'\in \mathcal{S}(\mathbb{R})$ and hence $\zeta(x,t)\in \mathcal{S}(\mathbb{R})$ for all $t$.

\section{Scaling of the Hamiltonian}

By the introduction of a small arbitrary constant parameter, $\delta$, defined by
\begin{equation}
\label{delta_def}
\delta=\frac{h_1}{\lambda}
\end{equation}
where $\lambda$ is the wavelength of the internal wave, meaning the scaled system will be considered as having \textit{long} waves, the variables will be scaled according to
\begin{align}
\eta\rightarrow\delta\eta, \quad
\xi\rightarrow \xi \quad \mbox{ and } \quad
\partial\rightarrow\delta\partial\notag
\end{align}
noting that the differential operators $\partial$ and $D$ are related by
\begin{equation}
D=-i\partial_x
\end{equation}
and the wave number $k:= 2\pi /\lambda$ is an eigenvalue of $D$ for the monochromatic linear waves in the form $e^{ikx}$ and therefore $\partial\sim\mathcal{O}(\delta)$.

The (unscaled) Dirichlet-Neumann operators can be expanded in terms of orders of $\eta$ as \cite{CraigSulem_NumSimofGravWaves}
\begin{align}
    G(\eta )&=|D|+D  \eta   D  -|D| \eta  |D|+\mathcal{O}(\eta ^2)\notag\\
\text{and} \quad G_1(\eta )&=D\tanh(h_1 D)-D  \eta  D  +D\tanh(h_1 D) \eta  D\tanh(h_1 D)+\mathcal{O}(\eta ^2)\notag
\end{align}
and so the expanded Dirichlet-Neumann operators are scaled as, noting from \cite{CraigGuyenneKalisch} that the constant term for the infinite lower layer is $|D|$,
\begin{align}
    G(\eta;\delta)&=\delta|D|+ \delta^3\big(  D   \eta D-   |D|  \eta   |D|\big)
    +\mathcal{O}( \delta^5) \notag\\
\text{and} \quad  G_1(\eta;\delta)&=\delta D \tanh(\delta h_1 D)\notag\\
  &\quad\quad\quad\quad-\delta^3\big(D\eta  D -  D \tanh(\delta h_1 D)  \eta  D \tanh(\delta h_1 D)\big)+\mathcal{O}( \delta^6). \notag
\end{align}

Using the expansion for the hyperbolic tangent it can be written that
\begin{equation}   
   \tanh(\delta h_1 D)=\delta h_1 D-\frac{1}{3}(\delta h_1 D)^3+\frac{2}{15}(\delta h_1 D)^5+\mathcal{O}(\delta^7),\notag
\end{equation}
and therefore the Dirichlet-Neumann operators can be expanded further as
\begin{align}
    G(\eta;\delta)&=\delta|D|+ \delta^3 D   \eta D- \delta^3  |D|  \eta   |D|
    +\mathcal{O}( \delta^5)\notag\\
\text{and} \quad G_1(\eta;\delta)&=\delta^2 h_1 D^2-\delta^3 D\eta D+\mathcal{O}( \delta^5).\notag
\end{align}

The details of the calculations are provided in the Appendix. The final expression for the Hamiltonian \eqref{Hamiltonian_conj}  is
\begin{multline}
\label{scaled_Ham}
H(\eta,\xi)=\frac{1}{2}\delta^2\frac{h_1}{\rho_1}\int\limits_{\mathbb{R}}\xi  D^2 \xi\,dx
-\frac{1}{2} \delta^3 \frac{\rho h_1^2}{\rho_1^2}\int\limits_{\mathbb{R}}\xi  |D|  D^2\xi\,dx
-\frac{1}{2}\delta^3\frac{1}{\rho_1}\int\limits_{\mathbb{R}}\xi  D\eta D\xi\,dx
\\
-\delta^3 \gamma_1\int\limits_{\mathbb{R}}\xi\eta\eta_x dx
-\delta^2 \kappa\int\limits_{\mathbb{R}} \xi\eta_x dx
+\frac{1}{6}\delta^3(\rho\gamma^2-\rho_1\gamma_1^2)\int\limits_{\mathbb{R}} \eta^3dx
+\frac{1}{2}\delta^2 A \int\limits_{\mathbb{R}}\eta^2 dx +\mathcal{O}( \delta^4)
\end{multline}
where 
\begin{equation}
A:=(\rho\gamma-\rho_1\gamma_1)\kappa+g(\rho-\rho_1 ).
\end{equation}

\section{The Benjamin-Ono approximation}

The variable
\begin{equation}
\mathfrak{u}=\xi_x\notag
\end{equation}
which assumes the role of momentum in the Hamiltonian approach (\textit{cf.} \cite{BenjBridPart1,BenjBridPart2}), in a similar fashion to $\eta$ assuming the role of the generalised coordinate, is introduced.

Noting that we can rescale the Hamiltonian by a factor $\delta^2$ by the choice of a proper time scale, the Hamiltonian (\ref{scaled_Ham}) in terms of $\eta$ and $\mathfrak{u}$ is given by
\begin{multline}
H(\eta, \mathfrak{u};\delta)=\frac{1}{2}\frac{h_1}{\rho_1}\int\limits_{\mathbb{R}}\mathfrak{u}^2\,dx
-\frac{1}{2} \delta \frac{\rho h_1^2}{\rho_1^2}\int\limits_{\mathbb{R}} \mathfrak{u}  |D|  \mathfrak{u}\,dx
-\frac{1}{2}\delta\frac{1}{\rho_1}\int\limits_{\mathbb{R}} \eta\mathfrak{u}^2\,dx
+\delta\frac{1}{2} \gamma_1\int\limits_{\mathbb{R}}\eta^2\mathfrak{u}\, dx\\
+ \kappa\int\limits_{\mathbb{R}}\eta\mathfrak{u} dx
+\frac{1}{6}\delta(\rho\gamma^2 - \rho_1\gamma_1^2) \int\limits_{\mathbb{R}} \eta^3dx
+\frac{1}{2} A \int\limits_{\mathbb{R}}\eta^2 dx +\mathcal{O}( \delta^2).
\end{multline}

The equations of motion \eqref{EOMsys} can be rewritten as
\begin{align}
\eta_t&=-\Big(\frac{\delta H}{\delta \mathfrak{u}}\Big)_x\notag \\
\text{and} \quad \mathfrak{u}_t&=-\Big(\frac{\delta H}{\delta \eta}\Big)_x-\Gamma\eta_t,\notag
\end{align}
and so therefore
\begin{align}
\label{BO_Eq1}
\eta_t&=- \frac{h_1}{\rho_1}\mathfrak{u}_x
+\delta \frac{\rho h_1^2}{\rho_1^2} |D|  \mathfrak{u}_x
+\delta\frac{1}{\rho_1}(\eta\mathfrak{u})_x
-\delta \gamma_1\eta\eta_x
- \kappa\eta_x \\
\label{BO_Eq2}
\text{and} \quad \mathfrak{u}_t&=\delta\frac{1}{\rho_1}\mathfrak{u}\mathfrak{u}_x
-\delta \gamma_1 (\eta\mathfrak{u})_x 
- \kappa\mathfrak{u}_x-\delta(\rho\gamma^2-\rho_1\gamma_1^2) \eta\eta_x
-A\eta_x -\Gamma\eta_t.
\end{align}

In the leading order \eqref{BO_Eq1} -- \eqref{BO_Eq2} give
\begin{align}
\eta_t+ \kappa\eta_x&=- \frac{h_1}{\rho_1}\mathfrak{u}_x
 \notag \\
\text{and} \quad \mathfrak{u}_t+ \kappa\mathfrak{u}_x&=
-A\eta_x -\Gamma\eta_t.\notag
\end{align}

This system of linear equations has a monochromatic solution (for a fixed wave number $k$) of the form 
\begin{align}
        \eta(x,t)&=\eta_0e^{ik(x-ct)}\notag\\
\text{and}   \quad     \mathfrak{u}(x,t)&=\mathfrak{u}_0e^{ik(x-ct)}\notag
\end{align}
where $c=c(k)$ is the wave speed. The linearised equations of motion produce
\begin{align}
\label{BO_c_Eq1}
-(c- \kappa) \eta&=- \frac{h_1}{\rho_1}\mathfrak{u}
\\
\label{BO_c_Eq2}
\text{and} \quad  -(c- \kappa) \mathfrak{u}&=
(-A+\Gamma c)\eta.
\end{align}

Multiplying both equations by each other gives their compatibility condition
\begin{align}
(c- \kappa)^2 =- \frac{h_1}{\rho_1}(-A+\Gamma c)\notag
\end{align}
which, in turn, gives the dispersion relation
\begin{align} \label{c}
c- \kappa=-\frac{h_1}{2\rho_1}\Gamma \pm\frac{1}{2}\sqrt{\frac{h_1^2}{\rho_1^2}\Gamma^2+4\frac{h_1}{\rho_1}g(\rho-\rho_1)}. \end{align}

As usual there are two solutions for the two signs $\pm$, which correspond to left and right-running waves with respect to the average flow with velocity $\kappa$.  Moreover, the leading order approximation from (\ref{BO_c_Eq1}) is
$$ \mathfrak{u}=\frac{\rho_1}{h_1}(c- \kappa) \eta$$ for a known $c$. Our aim now is to find the next order approximation of the form
\begin{equation} \label{u}
        \mathfrak{u}=\frac{\rho_1}{h_1}(c- \kappa) \eta+\delta\alpha\eta^2+\delta\beta |D|\eta + \mathcal{O}(\delta^2),
\end{equation}
for some, as yet, unknown values $\alpha$ and $\beta$, with a view to establishing a Benjamin-Ono type approximation.

Equation (\ref{BO_Eq1}) is hence written, to first order of $\delta$, as
\begin{multline}
\eta_t+ \kappa\eta_x=
 - \frac{h_1}{\rho_1}\bigg[\frac{\rho_1}{h_1}(c- \kappa) \eta+\delta\alpha\eta^2+\delta\beta |D|\eta\bigg]_x\\
+\delta \frac{\rho h_1^2}{\rho_1^2} |D|  \bigg[\frac{\rho_1}{h_1}(c- \kappa) \eta\bigg]_x
+\delta\frac{1}{\rho_1}\bigg[\eta\bigg(\frac{\rho_1}{h_1}(c- \kappa) \eta\bigg)\bigg]_x
-\delta \gamma_1\eta\eta_x
 \notag
\end{multline}
and so
\begin{equation}
\label{BO_alph_bet_Eq1}
\eta_t+c\eta_x +\delta\bigg[\frac{h_1}{\rho_1}\beta - \frac{\rho h_1}{\rho_1}(c- \kappa)\bigg]|D|\eta_x
+\delta\bigg[2\frac{h_1}{\rho_1}\alpha-2\frac{1}{h_1}(c- \kappa)+\gamma_1\bigg]\eta\eta_x=0.
\end{equation}

Likewise, equation (\ref{BO_Eq2}) is written as
\begin{multline}
\bigg[\frac{\rho_1}{h_1}(c- \kappa) \eta+\delta\alpha\eta^2+\delta\beta |D|\eta\bigg]_t+ \kappa\bigg[\frac{\rho_1}{h_1}(c- \kappa) \eta+\delta\alpha\eta^2+\delta\beta |D|\eta\bigg]_x+\Gamma\eta_t =\\
-A\eta_x+\delta\frac{1}{\rho_1}\bigg[\frac{\rho_1}{h_1}(c- \kappa) \eta\bigg]\bigg[\frac{\rho_1}{h_1}(c- \kappa) \eta\bigg]_x
-\delta \gamma_1 \bigg[\eta\bigg(\frac{\rho_1}{h_1}(c- \kappa) \eta\bigg)\bigg]_x 
-\delta(\rho\gamma^2-\rho_1\gamma_1^2) \eta\eta_x.
 \notag
\end{multline}

Noting that $\eta_t=-c\eta_x+\mathcal{O}( \delta)$
\begin{multline}
\label{BO_alph_bet_Eq2}
\eta_t +\bigg[\frac{\kappa \frac{\rho_1}{h_1}(c- \kappa)+A}{\frac{\rho_1}{h_1}(c- \kappa) +\Gamma}\bigg] \eta_x -\delta\bigg[\frac{\beta (c-\kappa)}{\frac{\rho_1}{h_1}(c- \kappa) +\Gamma}\bigg]|D| \eta_x \\
+\delta \bigg[ \frac{-2\alpha(c-\kappa)-\frac{\rho_1}{h_1^2}(c- \kappa)^2+2\gamma_1  \frac{\rho_1}{h_1}(c- \kappa)+\rho\gamma^2-\rho_1\gamma_1^2}{\frac{\rho_1}{h_1}(c- \kappa) +\Gamma}\bigg]\eta  \eta_x
=0.
\end{multline}

It is noted that the coefficient of the $\eta_x$ term 
\begin{equation}
\frac{\kappa \frac{\rho_1}{h_1}(c- \kappa)+A}{\frac{\rho_1}{h_1}(c- \kappa) +\Gamma} =c.
 \notag
\end{equation}

Comparing the $|D|\eta_x$ terms in (\ref{BO_alph_bet_Eq1}) and (\ref{BO_alph_bet_Eq2})
\begin{equation}
\delta\bigg[\frac{h_1}{\rho_1}\beta - \frac{\rho h_1}{\rho_1}(c- \kappa)\bigg]
 =
 -\delta\bigg[\frac{\beta (c-\kappa)}{\frac{\rho_1}{h_1}(c- \kappa) +\Gamma}\bigg]
 \end{equation}
and so
\begin{equation} \label{beta}
\beta
 =\frac{\rho\rho_1(c- \kappa)^2+\rho  h_1 \Gamma (c-\kappa)}{2\rho_1(c-\kappa)+h_1 \Gamma}.
\end{equation}

Comparing the $\eta\eta_x$ terms in (\ref{BO_alph_bet_Eq1}) and (\ref{BO_alph_bet_Eq2})
\begin{multline}
2\frac{ h_1}{\rho_1}\alpha-\frac{2}{h_1}(c- \kappa)+\gamma_1
=
\frac{-2\alpha(c-\kappa)-\frac{\rho_1}{h_1^2}(c- \kappa)^2+2 \frac{\rho_1\gamma_1 }{h_1}(c- \kappa)+\rho\gamma^2-\rho_1\gamma_1^2}{\frac{\rho_1}{h_1}(c- \kappa) +\Gamma}\notag
\end{multline}
and so
\begin{equation} \label{alpha}
\alpha
=
\rho_1\bigg(\frac{ \rho_1 (c- \kappa)^2
+2 h_1 \Gamma (c- \kappa)
- \gamma_1 h_1^2 \Gamma
+\rho_1 \gamma_1 h_1 (c-\kappa)
+h_1^2(\rho\gamma^2-\rho_1\gamma_1^2)}{2h_1^2 \big(2\rho_1(c- \kappa)+h_1\Gamma  \big)}\bigg).
\end{equation}

The equation for $\eta$, from (\ref{BO_alph_bet_Eq1}), is therefore given by
\begin{multline}
\eta_t+c\eta_x +\delta\Bigg[\frac{h_1}{\rho_1}\left(\frac{\rho\rho_1(c- \kappa)^2+\rho  h_1 \Gamma (c-\kappa)}{2\rho_1(c-\kappa)+h_1 \Gamma} \right) - \frac{\rho h_1}{\rho_1}(c- \kappa)\Bigg]|D|\eta_x\\
+\delta\Bigg[\rho_1\bigg(\frac{ \rho_1 (c- \kappa)^2
+2 h_1 \Gamma (c- \kappa)
- \gamma_1 h_1^2 \Gamma
+\rho_1 \gamma_1 h_1 (c-\kappa)
+h_1^2(\rho\gamma^2-\rho_1\gamma_1^2)}{2h_1^2 \big(2\rho_1(c- \kappa)+h_1\Gamma  \big)}\bigg)\\-2\frac{1}{h_1}(c- \kappa)+\gamma_1\Bigg]\eta\eta_x=0
\end{multline}
which can be written as
\begin{multline} \label{BO}
\eta_t+c\eta_x -\delta\Bigg[\frac{\rho h_1(c-\kappa)^2}{2\rho_1(c-\kappa)+h_1 \Gamma }\Bigg]|D|\eta_x\\
+\delta
\Bigg[\frac{-3\rho_1(c-\kappa)^2
+3\rho_1 \gamma_1 h_1(c-\kappa)
+ h_1^2 (\rho\gamma^2-\rho_1\gamma_1^2)}{h_1(2\rho_1(c-\kappa)+h_1\Gamma)}\Bigg]\eta\eta_x=0.
\end{multline}

The second component $\mathfrak{u}$ can be expressed with $\eta$ by \eqref{u} where now $\alpha$ and $\beta$ are given by \eqref{alpha} and \eqref{beta}.

We should keep in mind that there are always two sets of equations for the left and right running waves corresponding to the different choices of the sign in \eqref{c}. 

The obtained equation \eqref{BO} is the well known Benjamin-Ono (BO) equation \cite{BO1,BO2} which is an integrable equation whose solutions can be obtained by the Inverse Scattering method, see \cite{FA,KM98} and the references therein.

{\bf Remark}: In the absence of a current ($\kappa=0$, $\gamma=\gamma_1=0$) the equation becomes
\begin{equation}
\eta_t+c\eta_x -\frac{1}{2}\delta\frac{\rho h_1 c}{\rho_1 }|D|\eta_x
-\frac{3}{2}\delta
\frac{c
}{h_1} \eta\eta_x=0,
\end{equation} where

\begin{align}
c=\pm \sqrt{\frac{h_1}{\rho_1}g(\rho-\rho_1)}. \notag
\end{align}

\section{Solitary wave solution}

The \emph{standard} form of the BO equation is
\begin{equation}
        \eta_T +4\eta \eta_X +|\partial_X|\eta_{X} =0
\end{equation}
for which the one-soliton solution is known, \begin{equation} \label{ETA}
\eta(X,T)=\frac{C_0}{ C_0^2(X-C_0 T-X_0)^2+1},
\end{equation} where $C_0$ and $X_0$ are constants. Hence, the equation

\begin{equation}
        \eta_T +4\eta \eta_X +\sigma |\partial_X|\eta_{X} =0, \qquad  \sigma=\pm 1
\end{equation} has a solution 

$$\eta(X,T)=\frac{\sigma C_0}{ C_0^2(X-X_0 + \sigma C_0 T)^2+1},$$ which is the same solution, if one replaces the arbitrary constant $\sigma C_0 $ with $C_0$.

Transforming $X$ using $X\rightarrow X-cT $ gives the equation 
\begin{equation}
        \eta_T +c\eta_X+ 4\eta \eta_X +\sigma |\partial_X|\eta_{X} =0
\end{equation} with a solution 

$$\eta(X,T)=\frac{C_0}{ C_0^2[X-(c+C_0 )T-X_0]^2+1},$$ and hence after a rescaling of $\eta$, $X$ and $T$, equation \eqref{BO} has a solution

$$ \eta(x,t)=\frac{\eta_0}{C_0^2 \mu^2[x-x_0 -(c+C_0)t]^2+1} $$ where the amplitude $\eta_0$ and the initial displacement $x_0$ are arbitrary constants,  
\begin{equation} \label{C0} 
C_0:= \eta_0 \frac{-3\rho_1(c-\kappa)^2
+3\rho_1 \gamma_1 h_1(c-\kappa)
+ h_1^2 (\rho\gamma^2-\rho_1\gamma_1^2)}{4h_1(2\rho_1(c-\kappa)+h_1\Gamma)} \end{equation} 
and
\begin{equation}  
    \mu:= \frac{2\rho_1(c-\kappa)+h_1 \Gamma }{\rho h_1(c-\kappa)^2}.  \end{equation}

Expression \eqref{C0} shows that the wavespeed of the soliton depends on its amplitude $\eta_0.$  Also, note that the sign of $C_0$ depends on $\eta_0$ and the parameters of the system.

\section{Discussion}

The illustrative one-soliton solution of the BO equation \eqref{ETA} suffers, however, from the following disadvantages.  First, it is not in the Schwartz class in the $x$-variable, which is not a very serious disadvantage from the physical point of view. Second, it violates the assumption  \eqref{int_eta} for $\eta$ since for the expression \eqref{ETA}
$$ \int_{\mathbb{R}} \eta(X,T) dX = \pi \ne 0. $$ This is due to the fact that the initial condition $\eta(x,0)$ for this particular solution does not satisfy the mentioned assumptions.
Therefore, extra care is needed when the inverse scattering, or other methods are applied to the BO equation when modelling internal waves.

 
\section{Appendix}

Writing the Dirichlet-Neumann operators as
\begin{equation}
G(\eta;\delta)=\sum_{j=0}^{\infty} G^{[j]}(\eta;\delta)\mbox{ and }
G_1(\eta;\delta)=\sum_{j=0}^{\infty} G_1^{[j]}(\eta;\delta)
\end{equation}
where the superscript identifies the order of $\delta$ the relevant terms that we will be using are:
\begin{align}
    G^{[1]}&=\delta|D|,\quad G^{[2]}=0,\quad G^{[3]}=\delta^3 D   \eta D- \delta^3  |D|  \eta   |D| \notag\\
    \label{G_terms}
    G_1^{[1]}&=0, \quad G_1^{[2]}=\delta^2 h_1 D^2\mbox{ and } G_1^{[3]}=-\delta^3 D\eta D. 
\end{align}

Therefore noting that
\begin{align}
    \rho_1 G(\eta;\delta)&=\delta\rho_1 |D|+ \delta^3\rho_1 D   \eta D- \delta^3\rho_1  |D|  \eta   |D|+\mathcal{O}( \delta^5)\notag\\
\text{and} \quad \rho G_1(\eta;\delta)&=\delta^2\rho h_1 D^2-\delta^3 \rho D\eta D+\mathcal{O}( \delta^5)\notag
\end{align}
the operator $B$ is transformed to
\begin{equation}
B= \delta\rho_1 |D|+\delta^2\rho h_1 D^2+\mathcal{O}( \delta^3).\notag
\end{equation}

This can be written as
\begin{equation}
{B}=\delta \rho_1\frac{|D|}{D^2}D\Bigg\{1+\delta\frac{\rho}{\rho_1} h_1\frac{D^2}{|D|}+\mathcal{O}( \delta^2)\Bigg\}D\notag
\end{equation}
and so 
\begin{equation}
{B}^{-1}=\frac{1}{\delta\rho_1}\frac{D^2}{|D|}D^{-1}\Bigg\{1+\delta\frac{\rho}{\rho_1} h_1\frac{D^2}{|D|}+\mathcal{O}( \delta^2)\Bigg\}^{-1}D^{-1}.\notag
\end{equation}
Noting that
\begin{equation}
\frac{D^2}{|D|}=|D|\notag
\end{equation}
we can write
\begin{equation}
{B}^{-1}=\frac{1}{\delta\rho_1}|D|D^{-1}\Bigg\{1+\delta\frac{\rho}{\rho_1} h_1|D|+\mathcal{O}( \delta^2)\Bigg\}^{-1}D^{-1}.\notag
\end{equation}

Using the expansion
\begin{equation}
(1+x)^{-1}=1-x+\mathcal{O}(x^2)\notag
\end{equation}
the inverse of the operator $B$ is given by 
\begin{equation}
{B}^{-1}=\frac{1}{\delta\rho_1}|D|D^{-1}\bigg\{1-\delta\frac{\rho}{\rho_1} h_1|D|+\mathcal{O}( \delta^2)\bigg\}D^{-1}.\notag
\end{equation}

Writing the operator as
\begin{equation}
B^{-1}(\eta;\delta)=\sum_{j=0}^{\infty} B^{-1[j]}(\eta;\delta)
\end{equation}
where, again, the superscript identifies the order of $\delta$ the relevant terms that we will be using are:
\begin{equation}
    \label{B_terms}
    B^{-1[-1]}=\frac{1}{\delta\rho_1}|D|^{-1}\quad \mbox{and} \quad B^{-1[0]}=-\frac{\rho}{\rho_1^2} h_1.
\end{equation}

The Hamiltonian can therefore be written, using components of the expanded operators as
\begin{multline}\label{H-aux}
H(\eta,\xi)=\frac{1}{2}\int\limits_{\mathbb{R}}\xi G^{[1]}B^{-1[-1]}G_1^{[2]}\xi\,dx
+\frac{1}{2}\int\limits_{\mathbb{R}}\xi   G^{[1]}   B^{-1[0]}  G_1^{[2]} \xi\,dx\\
+\frac{1}{2}\int\limits_{\mathbb{R}}\xi   G^{[1]}   B^{-1[-1]}  G_1^{[3]} \xi\,dx
-\delta^3 \gamma\int\limits_{\mathbb{R}}\xi\eta\eta_x dx
-\delta^2 \kappa\int\limits_{\mathbb{R}} \xi\eta_x dx\\
+\delta^3\rho_1(\gamma-\gamma_1)\int\limits_{\mathbb{R}}\eta\eta_x  B^{-1[-1]}   G^{[1]} \xi\,dx
+\frac{1}{6}\delta^3(\rho\gamma^2-\rho_1\gamma_1^2)\int\limits_{\mathbb{R}} \eta^3dx\\
+\frac{1}{2}\delta^2 [(\rho\gamma-\rho_1\gamma_1)\kappa+g(\rho-\rho_1 )]\int\limits_{\mathbb{R}}\eta^2 dx +\mathcal{O}( \delta^4).
\end{multline}

Replacing the expressions for $G^{[j]}$, $G_1^{[j]}$ and $B^{-1[j]}$ and using (\ref{G_terms}) and (\ref{B_terms}) in \eqref{H-aux} gives the Hamiltonian \eqref{scaled_Ham}.


\section*{Acknowledgments} The authors would like to thank the Isaac Newton Institute for Mathematical Sciences, Cambridge, UK, for support and hospitality during the programme {\it Nonlinear water waves} where work on this paper was undertaken. The authors are thankful to two anonymous referees whose comments and suggestions have improved the quality of the manuscript.

This work was supported by EPSRC grant no EP/K032208/1. ACC is also funded by SFI grant 13/CDA/2117.

%

\medskip
Submitted to DCDS  August 2018; revised  October 2018.
\medskip

\end{document}